\input amstex.tex
\input amsppt.sty   
\magnification 1200
\vsize = 9.5 true in
\hsize=6.2 true in
\NoRunningHeads        
\parskip=\medskipamount
        \lineskip=2pt\baselineskip=18pt\lineskiplimit=0pt
       
        \TagsOnRight
        \NoBlackBoxes

        \topmatter
        \title
        Supercritical Nonlinear Schr\"odinger \\
        Equations II: Almost Global Existence 
        \endtitle
\author
         W.-M.~Wang        \endauthor        
\address
{D\'epartement de Math\'ematique, Universit\'e Paris Sud, 91405 Orsay Cedex, FRANCE}
\endaddress
        \email
{wei-min.wang\@math.u-psud.fr}
\endemail
\abstract
We prove almost global existence for supercritical nonlinear Schr\"odinger equations on the $d$-torus ($d$ arbitrary) on the good geometry selected in part I. This is seen as the Cauchy consequence of I,  since the known invariant measure of smooth solutions are supported on KAM tori. In the high frequency limit, these quantitative solutions could also be relevant to Cauchy problems for compressible Euler equations.

\endabstract

        \bigskip\bigskip
        \bigskip
        \toc
        \bigskip
        \bigskip 
        \widestnumber\head {Table of Contents}
        \head 1. Introduction and statement of the Theorem
        \endhead
        \head 2. Approximate solutions
        \endhead
        \head 3. Linearized flow
        \endhead
        \head 4. Almost global existence
        \endhead
        \endtoc
        \endtopmatter
        \vfill\eject
        \bigskip
\document
\head{\bf 1. Introduction and statement of the Theorem}\endhead
In this second part of the paper, we consider the Cauchy problems for the nonlinear Schr\"odinger 
equations (NLS) on the $d$-torus, $\Bbb T^d=[0, 2\pi)^d$. In accordance to the custom, we set the higher order 
terms $H$ in part I (I, 1.1) to be zero. This is also because
for small data, the construction below carries over verbatim to $H\neq 0$.

It is convenient to add a parameter and consider initial data of size one. So we have the following
Cauchy problem on $\Bbb T^d$:
$$
\cases i\frac\partial{\partial t}u =-\Delta u+\delta |u|^{2p}u\qquad (p\geq 1, p\in\Bbb N \text{ \it arbitrary}),\\
u(t=0)=u_0,\endcases
\tag 1.1
$$
with periodic boundary conditions: $u(t,x)=u(t, x+2n\pi)$, $x\in [0, 2\pi)^d$ for all $n\in\Bbb Z^d$
and $\delta\neq 0$ is the parameter.

It is well known, cf. \cite{Bo1} that (1.1) is locally wellposed in $H^s$ for 
$$s>\max (0, \frac{1}{2}(d-\frac{2}{p}))\tag 1.2$$
for a time interval of size $|\delta|^{-1}$.  When $d=p=1$, (1.1) is integrable and is locally and hence globally 
well-posed in $\Bbb L^2$. (1.2) is derived by linearizing about the flow of the Laplacian and using 
$L^p$ estimates of its eigenfunction solutions (Strichartz estimates). 

For $d\geq 3$ and sufficiently large $p$, the right side of (1.2) is larger than $1$. For example, in
dimension $3$, the quintic NLS is locally well-posed in $H^s$ for $s>1$, where there is no
conservation law. These equations are therefore {\it supercritical} as there is no a priori global 
existence from patching up local solutions, not even for small $\delta$, as (1.1) is non-dispersive,
i. e., $\Vert u\Vert_\infty$ cannot tend to $0$ as $t\to\infty$ on the torus $\Bbb T^d$. 

In part I, we went one step further and analyzed the resonance geometry created by the nonlinear
term $|u|^{2p}u$. Relying on the geometric information afforded by this analysis and linearizing
about a suitable approximate quasi-periodic solutions, we prove the following:
\proclaim{Theorem}
Let $u_0=u_1+u_2$. Assume $u_1$ is {\it generic} satisfying (I. i-iv.) and $\Vert u_2\Vert=\Cal O(\delta)$, 
where $\Vert\cdot\Vert$ is an analytic norm (about a strip of width $\Cal O(1)$) on $\Bbb T^d$. Let $\Cal B (0,1)=(0, 1]^b$,
where $b$ is the dimension of the Fourier support of $u_1$. Then there exists an open set $\Cal A\subset
\Cal B (0,1)$ of positive measure, such that for all $A>1$, there exists $\delta_0>0$, such that for 
all $\delta\in(-\delta_0, \delta_0)$, if $\{|\hat u_1|\}\in\Cal A$, then (1.1) has a unique solution $u(t)$ for 
$|t|\leq\delta^{-A}$ satisfying $u(t=0)=u_0$ and $\Vert u(t)\Vert  \leq \Vert u_0\Vert+\Cal O(\delta)$, moreover
meas $\Cal A\to 1$ as $\delta\to 0$.   
\endproclaim

\noindent{\it Remark.} It is essential that the set $\Cal A$ is {\it open}, as we will 
need to establish an open mapping theorem to analyze Cauchy problems .

For perturbations of the $1d$ cubic NLS ($d=p=1$), similar stability results are proven in \cite{Ba, Bo3}. For parameter
dependent equations see \cite{BG, Bo2}. The equations treated in \cite{Ba, BG, Bo2, 3} are either $\Bbb L^2$ 
or  essentially $\Bbb L^2$ well-posed. So there is a priori global existence. 

The equations treated in the Theorem are of a different nature, there is no a priori global existence
from conservation laws. In fact existence is obtained via explicit construction. As mentioned earlier
in I, this is possible because the known invariant measure for smooth flow is supported on 
KAM tori. Linearizing about approximate quasi-periodic solutions to prove existence and uniqueness for a 
time arbitrarily longer than local existence time is the main novelty of the paper. 

As in part I, we also have the following semi-classical counterpart, providing quantitative almost global $\Bbb L^2$ size $1$
and large (kinetic) energy solutions to Cauchy problems. These solutions could be relevant to Cauchy problems
for compressible Euler equations in the low density limit, cf. \cite{G, Se}.

\proclaim{Corollary}
Set $\delta=1$ in (1.1). Assume $u_0$ is {\it generic} with frequencies $\{j_k\}_{k=1}^b\in [K\Bbb Z^d]^b$, $K\in\Bbb N^+$.
Let $\Cal B (0,1)=(0, 1]^b$. Then there exists an open set $\Cal A\subset
\Cal B (0,1)$ of positive measure, such that for all $A>1$, there exists $K_0>0$, such that for 
all $K>K_0$, if $\{|\hat u_0|\}\in\Cal A$, then (1.1) has a unique solution $u(t)$ for 
$|t|\leq K^{-A}$ satisfying $u(t=0)=u_0$ and $\Vert u(t)\Vert  \leq \Vert u_0\Vert+\Cal O(1/K^2)$, 
where $\Vert\cdot\Vert$ is an analytic norm (about a strip of width $\Cal O(1/K)$) on $\Bbb T^d$, 
moreover meas $\Cal A\to 1$ as $K\to \infty$.   
\endproclaim

\demo{Proof}
Using the observations in the proof of the Corollary in part I, the proof follows that of the Theorem verbatim.
\hfill$\square$

\enddemo

\bigskip
\noindent{\it A sketch of the proof}

Writing the first equation in (1.1) as $F(u)=0$, for $u_0$ satisfying the conditions in the Theorem,
we first find an approximate solution $v$ such that 
$$\cases F(v)=\Cal O(\delta^r),  \qquad\qquad\qquad\qquad(1.3)\\
v(t=0)-u_0=\Cal O(\delta^r), \,\quad\quad\qquad\,(1.4)\endcases$$
where $r>A>1$.

This approximate solution $v$ is quasi-periodic with $\Cal O(|\log\delta|)$ number of basic frequencies.
Moreover at $t=0$, $v$ has the decomposition:
$$v(t=0)=u_1+v_2$$
with $u_1$ as in the Theorem and $\Vert v_2\Vert=\Cal O(\delta)$.

The construction of $v$ comprises of two steps. The first step is to construct approximate quasi-periodic
solutions of $\Cal O(|\log\delta|)$ number of basic frequencies with the initial approximation the
solution to the linear equation
$$u^{(0)}=u_1+u_2\tag 1.5$$
for all $u^{(0)}$ such that $u_1$ is generic and $\Vert u_2\Vert=\Cal O(\delta)$. 
This is done in section 2, using a finitely iterated Newton scheme. 

The excision in $\{|\hat u_1|\}$ is essentially the same as in sects. 2 and 3 of I, ensuring the existence of a spectral gap.
The amplitudes $\{|\hat u_2|\}$ are arbitrary as long as $\Vert u_2\Vert=\Cal O(\delta)$ is smaller
than the spectral gap. The constructed solution $u$ satisfies 
$$F(u)=\Cal O(\delta^r).\tag 1.6$$

Since the above construction is valid on a open set, using the open mapping theorem, we show in section 4
that for all 
$$\tilde u_0=u_1+\tilde u_2$$ 
of $\Cal O(|\log\delta|)$ number of frequencies, there is 
$$u^{(0)}=u_1+v_2$$ 
of the same number of frequencies such that the corresponding quasi-periodic solution $v$
satisfies 
$$\cases F(v)=\Cal O(\delta^r),  \\
v(t=0)=\tilde u_0.\endcases\tag 1.7$$

In section 3, we differentiate (1.6) with respect to the Fourier coefficients of $u^{(0)}$ in (1.5) and 
prove that the solutions to the linearized equation is a basis which spans $\Bbb L^2(\Bbb T^d)$,
after a further excision of $\{|\hat u_1|\}$. 
Schematically this could be understood as follows. 

Assume $u$ is a solution satisfying the equation 
$F(u)=0$ and that it depends on a parameter $a$, then $\partial u/\partial a$ is a solution 
to the linearized equation: $$F'(u)(\frac{\partial u}{\partial a})=0.$$
The fact that $\{\frac{\partial u}{\partial a}\}$ is a basis is a direct consequence of the separation property of
the resonance geometry entailed by generic $u_1$. This basis in turn allows us to control the flow
linearized about the $v$ in (1.7). Using Duhamel's formula and the linearized flow to control
the difference of (1.1) and (1.7), we conclude the proof of the Theorem.
\bigskip
\head{\bf 2. Approximate solutions}\endhead
We construct approximate solutions to  
$$ i\frac\partial{\partial t}u =-\Delta u+\delta |u|^{2p}u\qquad (p\geq 1, p\in\Bbb N)\tag 2.1$$
on the $d$-torus $\Bbb T^d$, which will serve as a ``basis" for the analysis of Cauchy problems in sect. 4.
Let $u$ be an entire function on $\Bbb T^d$: 
$$u(x)=\sum_{k\in G\subset\Bbb Z^d, \, |G|=b} a_k e^{-i\theta_k}e^{ij_k\cdot x},\tag 2.2$$
where $a_k\in\Bbb R^{*}_+$, $\theta_k\in[0, 2\pi)$ and $j_k\in\Bbb Z^d$. 

Let 
$$u(t, x)=\sum_{k\in G} a_k e^{-ij_k^2t}e^{-i\theta_k}e^{ij_k\cdot x}$$
be the space-time counter part of (2.2) which is a solution to (2.1) when $\delta=0$.
In view of the factor $\delta$ in (2.1),
we may assume without loss that $a_k\in (0, 1]$, so 
$$a=\{a_k\}_{k=1}^b\in (0, 1]^b=\Cal B(0, 1)=\Cal B,$$
where we have identified $G$ with the set $\{1, 2,...,b\}$.

Define $\Bbb A$ to be the class of analytic functions on $\Bbb T^d$ such that 
if $u\in\Bbb A$, then
\item {($\Bbb A1$)} $u$ has the decomposition: 
$$u=u_1+u_2,$$
where $u_1(\cdot)$ is entire (of form (2.2)) and $u_1(\cdot,\cdot)$ {\it generic} satisfying (I, i-iv) and 
$$a=\{a_k\}_{k=1}^b\in \Cal B\backslash\Cal B_\epsilon,$$
with $\Cal B_\epsilon$ as in Lemma 2.1 of part I, and

\item  {($\Bbb A2$)} $u_2=\Cal O(\delta)$ in an analytic norm:  $\exists \beta>0$, such that 
$$\sum_{j\in\Bbb Z^d} e^{\beta |j|} |\hat u_2(j)|<|\delta|.$$

Fix $\epsilon>0$, $0<\delta<\epsilon$ and $r>1$. Given $u\in\Bbb A$, let $B_J$ be the $\ell^\infty$
ball of radius $J$:
$$B_J=\{j\in\Bbb Z^d|\, \Vert j\Vert_\infty\leq J\},\tag 2.3$$
such that $\text{supp }u_1\subset B_J$ and 
$$u|_{B_J^c}=\Cal O(\delta^r)\tag 2.4$$
in the analytic norm ($\Bbb A2$), where 
$$u|_{B_J^c}:=\sum_{j\in B_J^c}\hat u(j) e^{ij\cdot x}. \tag 2.5$$
 
We make the identification:
$$\tilde u(t, x)|_{t=0}=u(x)\tag 2.6$$
by setting 
$$\hat { \tilde u}(-e_j, j)=\hat u(j)e^{-i\theta_j}
=a_j\in\Bbb R_+^{*}.\tag 2.7$$
Let $\tilde j\notin B_J$ and 
$$u^{(0)}(t, x)=\sum_{j\in B_J\cup\{\tilde j\}}a_je^{-i(\theta_j+j^2t)}e^{ij\cdot x},\tag 2.8$$
such that $$u^{(0)}(t, x)|_{t=0}\in\Bbb A.\tag 2.9$$
To economize, we will simply write $u^{(0)}\in\Bbb A$ when (2.8, 2.9) are verified.  

In the following, for all $\tilde j$, we construct approximate quasi-periodic solutions to (2.1) with $|B_J|+1=B$
frequencies. From (2.9), aside from a set $G\subset\Bbb Z^d$ of $b$ frequencies, $a_j=\Cal O(\delta)$
for $j\notin G$. In section 3, we eventually let $a_{\tilde j}\to 0$.

Taking $u^{(0)}$ as the initial approximation, we use as ansatz the nonlinear space-time Fourier series
$$u(t, x)=\sum_{(n,j) \in \Bbb Z^{B+d}}\hat u (n, j)e^{in\cdot (\theta+\omega t)}e^{ij\cdot x},\tag 2.10$$
where $\theta=\{\theta_j,\, j\in B_J\cup\{\tilde j\}\}$ and $\{\hat u(-e_j, j)\}=\{a_j\}$ are fixed;  while
$\omega=\{\omega_j,\, j\in B_J\cup\{\tilde j\}\}$ and $\{\hat u(n, j), \, (n, j)\neq (-e_j, j)\}$ are to be determined
during the construction.

(2.10) transforms (2.1) into a nonlinear matrix equation:
$$\text {diag }(n\cdot\omega+j^2)\hat u+\delta (\hat u*\hat {\bar u})^{*p}* \hat u=0.\tag 2.11$$
Let $\hat v=\hat{\bar u}$, writing also the equation for $\hat v$, we obtain as in part I:
$$
\cases
\text{diag }(n\cdot\omega+j^2)\hat u+\delta (\hat u*\hat v)^{*p}* \hat u=0,\qquad\qquad(2.12)\\
\text{diag }(-n\cdot\omega+j^2)\hat v+\delta (\hat u*\hat v)^{*p}* \hat v=0.\,\qquad\quad\,(2.13)
\endcases
$$
We write (2.12, 2.13) collectively as 
$$F(\hat u,\,  \hat v)=0,\tag 2.14$$
and use a modified Newton scheme to solve (2.14) to the desired order.

The main result of this section is the following:
\proclaim{Proposition 2.1}
Let $\Bbb A\ni u^{(0)}$ a solution to the linear equation of $B$ frequencies be as in (2.8). For
all $r\in\Bbb N^+$, there exists $\delta_0$. Let $a=\{a_k,\, k\in G\}$. There is $\Cal B'_\epsilon\supset \Cal B_\epsilon$,
meas $\Cal B'_\epsilon<\Cal O(\epsilon^c)$ ($c>0$), such that for all $a\in\Cal B\backslash \Cal B'_\epsilon$, all 
$\tilde j\notin B_J$, there is an approximate solution $u$ to (2.1) if $\delta$ satisfies $|\delta|<\delta_0$: 
$$i\dot{u}+\Delta u-\delta |u|^{2p}u=\xi,\tag 2.15$$
with $$u(t, x)=\tilde u(\theta+\omega t, x)=\tilde u(\tilde \theta, x),$$
$$\xi(t, x)=\tilde\xi (\theta+\omega t, x)=\tilde \xi(\tilde \theta, x).$$
$\tilde u$, $\tilde\xi$ are periodic on $\Bbb T^{B+d}$, entire in $\tilde\theta$ and analytic in $x$, satisfying
$$\Vert u-\sum_{j\in B_J\cup\{\tilde j\}}a_je^{-i(\theta_j+\omega_j t)}e^{ij\cdot x}\Vert\leq \Cal O(\delta)\tag 2.16$$
with $$\omega_j=j^2+\delta\Omega+\Cal O(\delta^2),\qquad j\in B_J\cup\{\tilde j\},\tag 2.17$$
and $\Omega$ depending only on $\{a_j, j\in G\}$ is the same for all $j$; 
$$\Vert\tilde\xi\Vert< \delta^r;\tag 2.18$$
$$\Vert \partial_a^\beta\partial_\omega^\gamma\tilde\xi\Vert <\delta^{r-\gamma};\tag 2.19$$
$$\partial_{\tilde\theta_j}\tilde\xi|_{a_j=0}=0\tag 2.20$$
for bounded $\beta$ and $\gamma$,  $a=\{a_j|j\in B_J\cup\{\tilde j\}\}$, $a$ and $\omega$ 
are viewed as independent variables in (2.19) and $\Vert\cdot\Vert$ is an appropriate analytic
norm (about a strip) in $t$ and $x$, uniformly in $\tilde j$.
\endproclaim

\noindent{\it Remark.} After differentiation with respect to $a_j$ and $\theta_j$, the $u$ constructed above provide a
 ``basis" for the linearized flow. So it is essential that they have localized Fourier support 
{\it uniformly} in $\tilde j$. We note also that Proposition 2.1 holds for {\it open} sets in $\{a_j\}_{j\in G}$ and
{\it all} $a_j$, $j\notin G$.  

\demo{Proof of Proposition 2.1}
For simplicity of notation, we drop the hat and write $u$ for $\hat u$ and $v$ for $\hat v$ until (2.41). Let 
$$\Cal S=\text{supp }u^{(0)}\cup \text{supp }v^{(0)}=\{(\mp e_j, j)\}_{j\in B_J\cup\{\tilde j\}}\subset\Bbb Z^{B+d}$$
be the resonant set. We use a Lyapunov-Schmidt decomposition to solve $F(u, v)=0$. $F$ restricted to $\Cal S$:
$F(u, v)|_{\Cal S}=0$ are the $Q$-equations;  $F(u, v)|_{\Bbb Z^{B+d}\backslash\Cal S}=0$ the $P$-equations.

Let $F'(u, v)$ be the linearized operator (evaluated at $u$, $v$):
$$\aligned 
F'(u, v) =&\pmatrix \text {diag }(n\cdot\omega+j^2)&0\\ 0& \text {diag }(-n\cdot\omega+j^2)\endpmatrix\\
&+\delta \pmatrix (p+1)(u*v)^{*p}& p(u*v)^{*p-1}*u*u\\ p(u*v)^{*p-1}*v*v& (p+1)(u*v)^{*p}\endpmatrix
\quad  (p\geq 1)\endaligned\tag 2.21$$
on $\ell^2(\Bbb Z^{B+d})$.  Given a set $\Bbb Z^{B+d}\supset A$, we define the restricted operator $F'_A$ to be:
$$\align F'_A(n, j; n', j')&=F'(n, j; n', j')\qquad \text {if } (n, j) \text{and } (n',j')\in A\tag 2.22\\
&=0 \qquad\qquad\qquad\quad  \text{ otherwise}.\tag 2.23\endalign$$
Let $$T=F'_{\Bbb Z^{B+d}\backslash\Cal S}\tag 2.24$$ 
and $T_N$ the restricted operator:
$$\aligned T_N(n, j;n', j')=&T(n, j;n', j')\qquad \text{if } |n|_1\leq N \text { and }|n'|_1\leq N;\\
=&0\qquad\qquad\qquad\quad\text{otherwise},\endaligned \tag 2.25$$
where $N=N(r)>r$. 

For $\Bbb A\ni u^{(0)}$, $u^{(0)}=u_1+u_2=\Cal O(1)+\Cal O(\delta)$,  
since the $\Cal O(1)$ component is generic, treating the $\Cal O(\delta)$ component as perturbation
we see that $T(u^{(0)})$ and $T_N(u^{(0)})$ are invertible with $\Vert T^{-1}(u^{(0)})\Vert$ and $\Vert T_N(u^{(0)})\Vert\leq\Cal O(\delta^{-1})$, where $\omega$ is set to 
$$\omega^{(0)}=\{j^2\}_{ j\in B_J\cup\{\tilde j\}}.\tag 2.26$$ 
We write $u$ collectively for $(u, v)$.
\bigskip
\noindent{\it Newton scheme}

To solve the $P$-equations, define
$$\Delta u^{(k)}=-T_N^{-1} (F(u^{(k-1)}))|_{\Bbb Z^{B+d}\backslash\Cal S}\tag 2.27$$
with the initial approximation:
$$\aligned u^{(0)}&=a_j\qquad \text {if } \pm j\in B_J\cup\{\tilde j\}\text{ and }n=\mp e_j\\
&=0\qquad\,\,\text{otherwise}.\endaligned\tag 2.28$$
and $$u^{(k)}=u^{(k-1)}+\Delta u^{(k)}\qquad k\geq 1,\tag 2.29 $$
the $k^{\text {th}}$ approximation. Below for notational simplicity, we often drop the restriction
sign: $|_{{\Bbb Z^{B+d}}\backslash\Cal S}$.

On $\Cal S$, $u$ is held fixed
$$u|_{\Cal S}=u^{(0)}\qquad\text{ as in } (2.28).\tag 2.30$$
We solve the $Q$-equations by solving for the $B$ frequencies:
$$\omega_j^{(k)}=j^2+\frac{\delta}{a_j} (u^{(k-1)}*v^{(k-1)})^{*p}*u^{(k-1)}|_{(-e_j, j)},\qquad j\in B_J\cup\{\tilde j\},\quad k\geq 1\tag 2.31$$
with the initial approximation, the eigenvalues of the linear operator:
$$\omega_j^{(0)}=j^2,\qquad j\in B_J\cup\{\tilde j\}.\tag 2.32$$
We first solve the $P$-equations according to (2.27) with $\omega=\omega^{(0)}$.

Since $u^{(0)}\in\Bbb A$, using (2.27) we obtain 
$$\Delta u^{(1)}=\Cal O(\delta).\tag 2.33$$
The remainder 
$$F(u^{(1)})=F(u^{(0)}+\Delta u^{(1)})=F(u^{(0)})+F'(u^{(0)})\Delta u^{(1)}+\delta\Cal O((\Delta u^{(1)})^2).$$
Since 
$$T_N^{-1}=T^{-1}+T^{-1} (T-T_N) T_N^{-1},\tag 2.34$$
we obtain from (2.34) 
$$F(u^{(1)})=\Cal O(\delta^3),\tag 2.35$$
where we also used exponential decay of matrix elements of 
$T^{-1}$, $T_N^{-1}$ and $(T-T_N)$.
Using $u^{(0)}$ in (2.31), we obtain 
$$\omega_j^{(1)}=j^2+\Cal O(\delta).\tag 2.36$$

Let $\omega=\{\omega_j^{(1)}\}$, $u=u^{(1)}=u^{(0)}+\Delta u^{(1)}$, in order to continue the iteration, we need
to estimate $T_N^{-1}(u^{(1)})$, where 
$$\aligned 
T_N(u^{(1)}) =&\pmatrix \text {diag }(n\cdot\omega+j^2)&0\\ 0& \text {diag }(-n\cdot\omega+j^2)\endpmatrix\\
&+\delta \pmatrix (p+1)(u*v)^{*p}& p(u*v)^{*p-1}*u*u\\ p(u*v)^{*p-1}*v*v& (p+1)(u*v)^{*p}\endpmatrix
\quad  (p\geq 1)\endaligned\tag 2.37$$
on $\ell^2(\{n, \, |n|_1\leq N\}\times\Bbb Z^d)$ with $\omega=\omega^{(1)}$, $u=u^{(1)}$ and $v=v^{(1)}$. We 
proceed as in section 3 of I, since $\omega_j^{(1)}=j^2+\Cal O(\delta)$ and $|n|_1\leq N$, the resonance structure
remains the same as in the first step when $\omega=\omega^{(0)}$ and $u=u^{(0)}$.  

From (2.31), we notice further that 
$$\omega_j^{(1)}=j^2+\delta\Omega+\Cal O(\delta^2),\qquad j\in B_J\cup\{\tilde j\},\tag 2.38$$
where $\Omega$ is the same for all $j$ and moreover only depends on $\{a_j|j\in G\}$. Treating 
the $\Cal O(\delta^2)$ terms in (2.38) and $\Delta u^{(1)}$ as perturbations, as in the proof of 
Lemma 3.1 of I, the invertibility of $T_N$ is therefore reduced to controlling the determinant of
matrices $\Gamma$ of sizes at most $(2|G|+d)\times (2|G|+d)$:
$$\Gamma=\pmatrix \text{diag }(n\cdot\Bbb I)\Omega &0\\0&\text{diag }(-n\cdot\Bbb I)\Omega \endpmatrix+\Cal A, \tag 2.39$$
where $\Cal A$ is a convolution matrix, which only depends on $\{a_j|j\in G\}$.

Since $\det\Gamma$ is a polynomial in $a_j$ of degree at most $2p(2|G|+d)$ and there are finite number of such
$\Gamma$'s depending only on $|G|$, $d$ and $N$ and independent of $\delta$, we have that there exist $C$, $c>0$,
such that for all $0<\epsilon<1$
$$\text{meas } \{a\in\Cal B|\,|\det\Gamma|<\epsilon, \, \text{ all }\Gamma\}\leq C\epsilon^c,\tag 2.40$$
where $a=\{a_j, j\in G\}$ and $\Cal B=(0, 1]^{|G|}=(0, 1]^b$. (This is the same as in the first part of the proof of Lemma 3.1 in I.) 

Let $\Cal B'_\epsilon$ be the set defined in (2.40), $\Cal B'_\epsilon\supset\Cal B_\epsilon$ as $\Cal B_\epsilon$ 
corresponds to taking $n=0$ in (2.39). We obtain that for $a\in\Cal B\backslash \Cal B'_\epsilon$,
$$\Vert T_N^{-1}(u^{(1)})\Vert =\Cal O(\delta^{-1}).\tag 2.41$$
Using (2.35, 2.27), iterating $r'$ times for some $r'<r$, putting back the hat and let $\hat u=u^{(r')}$ and $\hat\xi=F(u^{(r')})$,
we obtain (2.16-2.18). 

(2.20) follows from (2.10) by inspection. So we are only left to prove the derivative estimates in (2.19). We use induction. We have
$$\aligned &\partial_a F(u+\Delta u)\\
=&F'(u+\Delta u)\partial_a (u+\Delta u)\\
=&[F'(u)+F''(u)\Delta u+F'''(c)(\Delta u)^2][\partial_a u+\partial_a \Delta u]\\
=&\partial_a (F(u)+F'(u)\Delta u)\\
&+F''(u)\Delta u\partial_a\Delta u+F'''(c)(\Delta u)^2\partial_a u+F'''(c)(\Delta u)^2\partial_a \Delta u,\endaligned\tag 2.42$$
where we used Taylor and $c$ is between $u$ and $u+\Delta u$.

Since $$\Delta u=-T_N^{-1}(u)|_{\Cal S^c} F(u)|_{\Cal S^c},\tag 2.43$$
$$\aligned \partial_a\Delta u=&-T_N^{-1}(u)\partial_a F_N'(u)T_N^{-1}(u)|_{\Cal S^c}F(u)|_{\Cal S^c}\\
&-T_N^{-1}(u)|_{\Cal S^c}\partial_aF(u)|_{\Cal S^c}.\endaligned\tag 2.44$$
Using $\partial_a F_N'(u)=F_N''(u)\partial_a u$ and starting from $u=u^{(0)}$, we have 
$$\partial_a F_N'(u^{(0)})=\Cal O(\delta),\, \partial_a F(u^{(0)})|_{\Cal S^c}=\Cal O(\delta)\text{ and supp } \partial_a F(u^{(0)})
=\text{supp }F(u^{(0)}).\tag 2.45$$
So $$\partial_a\Delta u^{(1)}=\Cal O(\delta)\tag 2.46$$
and from (2.42) 
$$\partial_aF(u^{(1)})=\Cal O(\delta^3)=\Cal O(F(u^{(1)}).\tag 2.47$$
Iterating we obtain 
$$\Vert \partial_a\tilde\xi\Vert <\delta^r.\tag 2.48$$

To estimate the $\omega$ derivatives, we note that $F|_{\Cal S}=0$ identically and 
$$\frac{\partial}{\partial\omega}F(u^{(0)})|_{\Cal S^c}=0.\tag 2.49$$
From (2.43)
$$\aligned F(u+\Delta u)=&(F'-F_N')T_N^{-1}|_{\Cal S^c}F(u)|_{\Cal S^c}+F''(c)(\Delta u)^2\\
:=&A+F''(c)(\Delta u)^2,\endaligned\tag 2.50$$
for some $c$ between $u$ and $u+\Delta u$. 

So $\Cal F:=F(u+\Delta u)-A$ has no explicit dependence on $\omega$ and the derivatives can be estimated 
as in (2.42-2.47). Using 
$\partial_\omega F_N'(u^{(0)})=\Cal O(1)$ in the analogue of (2.44) to estimate  $\partial_\omega\Delta u$
and $\partial_\omega A=\Cal O(\delta^{N'})\partial_\omega F(u)$ ($N'>r$), we obtain 
$$\Vert \partial_\omega\tilde\xi\Vert <\delta^{r-1}.$$
Clearly higher order and mixed derivatives can be similarly estimated and we obtain (2.19).
\hfill $\square$
\enddemo
\bigskip
\head{\bf 3. Linearized flow}\endhead
The construction leading to the proof of Proposition 2.1 has a two fold application to Cauchy problems, since it is
valid on {\it open} sets of $a$. Firstly it is used to control solutions with initial conditions arbitrarily close to the given
one. We remark that from (2.16), the solutions in Proposition 2.1 are only $\Cal O(\delta)$ close to the given one. 
However, fundamentally due to an open mapping theorem that we establish in section 4, we will construct solutions that are $\Cal O(\delta^r)$ close to the given initial condition. Together with the estimate on the linearized flow, we then prove the Theorem.

\noindent{\it The linearized flow} 

From (2.15), we have approximate solutions $u$ to the nonlinear equation (2.1) satisfying 
$$F(u)=\xi.\tag 3.1$$
Differentiating with respect to $a$ and $\theta$:
$$\align&F'(u)\big(\frac{\partial u}{\partial a}\big)=\frac{\partial\xi}{\partial a}+\frac{\partial\xi}{\partial\omega}\cdot \frac{\partial\omega}{\partial a}+\frac{\partial\xi}{\partial\theta}\cdot \frac{\partial\omega}{\partial a}\cdot t,\qquad \{(a_j, \theta_j), j\in B_J\cup\{\tilde j\}\},\tag 3.2
\\
&F'(u)\big(\frac{\partial u}{\partial\theta}\big)=\frac{\partial\xi}{\partial\theta}, \qquad\qquad\qquad\qquad\qquad\qquad\quad \{(a_j, \theta_j), j\in B_J\cup\{\tilde j\}\},\tag 3.3\endalign$$
then letting $a_{\tilde j}\to 0$ and using (2.19, 2.20), we obtain that $\partial u/\partial a$ and $\partial u/\partial\theta$ are 
approximate solutions to the linear equation:
$$i\frac{\partial}{\partial t} w=-\Delta w+[\delta (p+1)|u(a_{\tilde j}=0)|^{2p}]w+[\delta p|u(a_{\tilde j}=0)|^{2(p-1)}u^2(a_{\tilde j}=0)]\bar w,\tag 3.4$$
where $u$ as in (2.15). (3.4) is the original nonlinear equation (2.1) linearized at $u(a_{\tilde j}=0)$.

\noindent{\it Remark.} It will become clear in the proof of Lemma 3.1 that (3.4) is indeed the corresponding linear equation due to continuity properties. 

We note that (3.4) is only linear over $\Bbb R$, so we need both  $\partial u/\partial a$ and $\partial u/\partial\theta$ to expand $\Bbb L^2$ solutions. Let  
$$\cases \nu _j(x, t)=\frac {1}{a_j} \frac{\partial u}{\partial\theta_j}(x, t)|_{a_{\tilde j}=0},\qquad\qquad\, (3.5)\\
 w_j(x, t)=\frac {\partial u} {\partial a_j} (x, t)|_{a_{\tilde j}=0},\qquad\qquad\quad (3.6)\endcases$$
where $u$ as in Proposition 2.1, $j\in B_J\cup \{\tilde j\}$, $B_J$ as in (2.3) and varying over all $\tilde j\notin B_J$.
Let $\nu_j^{(0)}=\nu_j(x, 0)$ and $w_j^{(0)}=w_j(x, 0)$. We have the following.
\proclaim{Lemma 3.1}
There is $\Cal B''_\epsilon\supset \Cal B'_\epsilon$,
meas $\Cal B''_\epsilon<\Cal O(\epsilon^c)$ ($c>0$), such that for all $a\in\Cal B\backslash \Cal B'_\epsilon$, the sequences $\{\nu_j^{(0)}\}$ and $\{w_j^{(0)}\}$ span $\Bbb L^2(\Bbb T^d)$. Moreover the flow $S(t)$ of (3.4) satisfies 
$$\Vert S(t)\Vert <1+2|t|\tag 3.7$$ for $t<\delta^{-r/3}$, where $\Vert\cdot\Vert$ is the analytic norm in $(\Bbb A2)$ with 
$\beta'\in (0, \beta)$ replacing $\beta$. 
\endproclaim
\demo{Proof}
From (2.27-2.29, 2.31) in the proof of Proposition 2.1, the Fourier coefficients $\hat u$ of the approximate solution $u$ to (2.1) is constructed using a Newton
scheme:
$$\align\Delta {\hat u}^{(k)}&=-T_N^{-1}({\hat u}^{(k-1)})F({\hat u}^{(k)}),\quad k\geq 1,\tag 3.8\\
 {\hat u}^{(k)}&={\hat u}^{(k-1)}+\Delta {\hat u}^{(k)},\\
 {\hat u}^{(0)}&=a_j, \quad\text{if } j\in B_J\cup\{\tilde j\}\text{ and }n=-e_j,\\
 &=0\quad\,\text{ otherwise},\\
 \omega_j^{(k)}&=j^2+\frac{\delta}{a_j} ({\hat u}^{(k-1)}*{\hat v}^{(k-1)})^{*p}*{\hat u}^{(k-1)}|_{(-e_j, j)},\qquad j\in B_J\cup\{\tilde j\},\quad k\geq 1\tag 3.9\\
 \hat u&={\hat u}^{(0)}+\sum_{k=1}^K\Delta {\hat u}^{(k)},\tag 3.10\endalign$$
 according to (2.10) using (3.8-3.10). Without loss we may assume in $(\Bbb A1)$ that 
 $$\text {supp }\hat u_1\cap \text{supp }\hat u_2=\emptyset.$$
 
We first prove  that $\{w_j^{(0)}\}$ spans $\Bbb L^2(\Bbb T^d)$. From (2.10), 
$$\aligned w_{j'}=\frac{\partial \hat u}{\partial a_{j'}}=&\sum\frac{\partial\hat u}{\partial a_{j'}}(n, j)e^{in\cdot(\theta+\omega t)e^{ij\cdot x}}\\
&+\sum\frac{\partial \hat u}{\partial \omega}\cdot \frac{\partial\omega}{\partial a_{j'}}e^{in\cdot(\theta+\omega t)e^{ij\cdot x}}\\
&+\sum\hat u(n, j)(in\cdot \frac{\partial\omega}{\partial a_{j'}} t)e^{in\cdot(\theta+\omega t)e^{ij\cdot x}}\\
&:=f_{j'}+\gamma_{j'}+g_{j'}.\endaligned\tag 3.11$$

\item{(i)} Estimates on $f$.

Let $\Cal C$ be the characteristic variety:
$$\{(n, j)|n\cdot\omega^{(0)}\pm j^2=0, \omega_j^{(0)}=j^2\}\tag 3.12$$ 
and $\Gamma=\{(\Delta n, \Delta j)\}$ the Fourier support of 
$|u_1^{(0)}|^{2p}$ and $|u_1^{(0)}|^{2(p-1)}{u_1^{(0)}}^2$ as in (I, 1.16-1.18), $u_1$ as in
$(\Bbb A1)$ is generic.

Let $\alpha\subset \Cal C$ ($|\alpha|\geq 2$), we say $\alpha$ is connected if for all $\alpha_1\in\alpha$,
$\exists\alpha_2\in\alpha$, $\alpha_1\neq\alpha_2$, such that $\alpha_1-\alpha_2\in\Gamma$.  Let $A'$
be the subset of connected sets on $\Cal C$ such that if $\alpha\in A'$ then there is an unique element in
$\alpha$ of the form $(\mp e_j, j)$. From sect. 2 of I, 
$$|\alpha|\leq 2b+d,\tag 3.13$$ 
where $b$ is the dimension of the Fourier support of $u_1$, as the notion of connectedness remains the same.

For $j'$ such that 
$$\nexists \alpha\in A' \text{ such that }(-e_{j'}, j')\in \alpha\tag 3.14$$
$$f_{j'}(t, x)=e^{-i(\theta_{j'}+\omega_{j'}t)}e^{ij'\cdot x}+\Cal O(\delta),\tag 3.15$$
where $\Cal O(\delta)$ is in the analytic norm in $(\Bbb A2)$ with the weight $e^{\beta'|j-j'|}$ replacing
$e^{\beta|j|}$ ($0<\beta'<\beta$) and is uniform in $t$.
Let $\text {supp }u_1$ be the Fourier support of $u_1$ and $G=\Pi_{\Bbb Z^d}\text {supp }u_1$ as before, where
$\Pi_{\Bbb Z^d}$ denotes projection onto $\Bbb Z^d$. we 
note that if $j\in G$, then $f_j$ satisfies (3.15). 

For $j'$ such that (3.14) is violated, let $\alpha\in A'$ be the maximal connected set containing $(-e_{j'}, j')$. 
Let $$\tilde\alpha=\Pi_{\Bbb Z^d}\alpha.\tag 3.16$$
From translation invariance for all $j\in\tilde\alpha$, there is (maximal) $\alpha(j)\subset\Cal C$ such that 
$(-e_j, j)\in\alpha(j)$ and $\Pi_{\Bbb Z^d}\alpha(j)=\tilde\alpha$.

Let $$M=F'(u_1^{(0)}, v_1^{(0)})|_{\alpha(j)},\quad j\in\tilde\alpha. \tag 3.17$$
$M$ is a $|\alpha|\times |\alpha|$ convolution matrix and hence the same for all $\alpha(j)$.
We have 
$$\aligned \frac{\partial\hat u}{\partial a_j}=&-M^{-1}|_{\alpha(j)\backslash (-e_j, j)}\delta (\hat u_1*\hat v_1)^{*p}\eta_{(-e_j, j)}|_{\alpha(j)\backslash (-e_j, j)}\\
&+\eta_{(-e_j, j)}+\Cal O(\delta),\quad\text{all }j\in\tilde\alpha\endaligned\tag 3.18$$
where $$\cases \eta_{(-e_j, j)}(n', j')=&1\quad (n',j')=(-e_j, j),\\
&0\quad \text{otherwise},\endcases$$
and $\Cal O(\delta)$ is in the analytic norm in $(\Bbb A2)$ (centered at $j$). 

We note that since $u_1^{(0)}$ is generic, $\Vert M^{-1}\Vert=\Cal O(\delta^{-1})$, 
moreover the matrix elements of $M^{-1}$ are rational functions of $\{a_j,\, j\in G\}$.
Let 
$$\align\phi_j&=e^{-i(\theta_j+\omega_j t)}e^{ij\cdot x}\\
&:=e^{-i\tilde \theta_j}e^{ij\cdot x},\quad j\in\tilde\alpha\tag 3.19\\
\phi&=\{\phi_j\}\\
f&=\{f_j\},\quad f_j \text{ as defined in }(3.11)\\
U&\text{ a unitary diagonal matrix:}\\
U&=\{e^{in'\cdot\tilde\theta}, \, (n', j') \in \alpha(j),\, j\in\tilde\alpha\}.\tag 3.20
\endalign$$
(The choice of $j$ is immaterial as it just produces an overall phase.)

Using (3.18), the above considerations give 
$$f=[U^{-1}\Cal M U]\phi+\Cal O(\delta),\tag 3.21$$ 
where $\Cal O(\delta)$ is in the analytic norm in $(\Bbb A2)$; $\Cal M_{jj}=1$ and $\Cal M_{jk}(j\neq k),\, j,\,k\in\tilde\alpha$
are rational functions in $\{a_j, j\in G\}$. 

\noindent{\it Remark.} The matrix $\Cal M$ is {\it independent} of $\theta$ due to genericity condition (i) in I, namely to 
this order $(\Delta n, 0)$ is not in the algebra generated by $\Gamma$ when restricted to $\Cal C$. 

We make an additional excision of the set 
$$M_\alpha=\{\{a_j,\, j\in G\}|\,|\det\Cal M|<\epsilon\},$$
with
$$\text{meas }M_\alpha<C'\epsilon^c\, (c>0).\tag 3.22$$
Since there are only finite numbers of different $\Cal M=\Cal M_\alpha$, we obtain from (3.22) that 
$$M=\cup_\alpha M_\alpha\tag 3.23$$
satisfies 
$$\text{meas }M_\alpha<C\epsilon^c\, (c>0).\tag 3.24$$
Set $\Cal B''_\epsilon=\Cal B'_\epsilon\cup M$, where $\Cal B_\epsilon$ as in Lemma 2.1 of I. For 
$a=\{a_j,\, j\in G\}\in\Cal B\backslash\Cal B''_\epsilon$,
$$\Vert \Cal M_\alpha^{-1}\Vert <\epsilon^{-1}\tag 3.25$$
for all $\alpha$. 

\item{(ii)} Estimates on $g$

We prove $$\frac{\partial\omega_k}{\partial a_j'}=\delta\Cal O(a_{j'})\tag 3.26$$
by induction. We have 
$$\omega_k^{(1)}=k^2+\frac{\delta}{a_k} ({\hat u}^{(0)}*{\hat v}^{(0)})^{*p}*{\hat u}^{(0)}|_{(-e_k, k)}.\tag 3.27$$

In view of the restriction to $(-e_k, k)$, the second term is even in $a_{j'}$, so 
$$\frac{\partial\omega_k^{(1)}}{\partial a_j'}=\delta\Cal O(a_{j'}).\tag 3.28$$
Assume $$\frac{\partial\omega_k^{(m-1)}}{\partial a_j'}=\Cal O(a_{j'}),$$
then at the $m^{\text {th}}$ iteration
$$\omega_k^{(m)}=k^2+\frac{\delta}{a_k} ({\hat u}^{(m-1)}*{\hat v}^{(m-1)})^{*p}*{\hat u}^{(m-1)}|_{(-e_k, k)}.\tag 3.29$$
The second term is even in $a_{j'}$, using (3.28) on the term which does not have explicit dependence on $a_{j'}$, we
obtain (3.26).

Using (3.26) in (3.11), we obtain 
$$g_{j'}=\delta\Cal O(a_{j'})\cdot t \Cal O(\sum|\hat u(n, j)|),\tag 3.30$$
where we used $\Vert n\Vert_{\ell ^1}\leq N$. Since 
$$|a_{j'}|\leq e^{-\beta |j'|}\tag 3.31$$ 
and from Proposition 2.1
$$\sum |\hat u(n, j)|e^{\beta' |j'|}<\infty,\tag 3.32$$
for some $0<\beta'<\beta$, there exists 
$0<\beta''<\beta-\beta'$ such that 
$$|g_{j'}e^{\beta''|j'|}|=\Cal O(\delta)\cdot t\tag 3.33$$ 
uniformly in $j'$. Similarly we estimate $\gamma_{j'}$ to be of $\Cal O(\delta)$:
$$\gamma_{j'}=\Cal O(\delta).\tag 3.34$$

So $w_{j'}(t)$ satisfies 
$$\sup_{j'}\sum_{|j-j'|>A}e^{\beta''|j-j'|}|\hat w_{j'}(j)|=\Cal O(\delta)\cdot t\tag 3.35$$
for some $A$ which only depends on the Fourier support of $u_1$. Using (3.15, 3.21, 3.33-3.35)
in (3.11), we conclude that $\{w_j^{(0)}\}$ spans $\Bbb L^2(\Bbb T^d)$. 

Similar arguments together with the observation that 
$$\frac{\partial u}{\partial \theta_{j'}}|_{a_j'=0}=0$$
prove that $\{\nu_j^{(0)}\}$ spans $\Bbb L^2(\Bbb T^d)$,
moreover $$\nu_{j}^{(0)}=iw_j^{(0)}+\Cal O(\delta),\tag 3.36$$
where $\Cal O(\delta)$ is in the analytic norm in $(\Bbb A2)$ with $e^{\beta''|j-j'|}$ replacing 
$e^{\beta|j|}$ ($0<\beta''<\beta$).

To prove (3.7), we first estimate the error terms in (3.2, 3.3). Let $\Phi$ denote either $\nu_j$ or $w_j$. It follows from (2.18-2.20, 3.26)
that 
$$i\Phi_t+\Delta \Phi-[\delta (p+1)|u|^{2p}]\Phi-[\delta p|u|^{2(p-1)}u^2]\bar \Phi=R,\tag 3.37$$
where $u$ is the approximate solution constructed in Proposition 2.1 with $a_{\tilde j}=0$ and 
$$\Vert R\Vert\leq (1+|t|)\delta^r\tag 3.38$$
uniformly in $j$, where $\Vert\cdot\Vert$ is the analytic norm defined in $(\Bbb A2)$ with $\beta'\in(0,\beta)$ replacing $\beta$.  
Here we also used $$\frac{\partial\xi}{\partial a_{\tilde j}}=\Cal O(a_{\tilde j})=0\text{ at } a_{\tilde j}=0.$$
This is because $\text{supp }\xi\subseteq\Bbb Z^{B+d}\backslash\Cal S$, since on $\Cal S$, $\xi=0$ identically, so
$n_{\tilde j}\neq \pm 1$ on $\text{supp }\xi$ and $\xi=\Cal O(a_{\tilde j}^2)$. 

Since $\{\nu_j^{(0)}\}$, $\{w_j^{(0)}\}$ span $\Bbb L^2(\Bbb T^d)$ and the left side of (3.37) is linear over $\Bbb R$, given any
$\Psi$ with $\Vert \Psi\Vert\leq 1$, there is the expansion 
$$\psi=\psi(x)=\sum \alpha_j\nu_j^{(0)}+\beta_jw_j^{(0)}\tag 3.39$$
with $\alpha_j$, $\beta_j\in\Bbb R$ and 
$$\Psi(t)=\sum \alpha_j\nu_j(t)+\beta_jw_j(t)\tag 3.40$$
satisfies (3.37. 3.38) and
$$\aligned \Vert \Psi(t)\Vert&\leq (1+\delta |t|)\Vert\Psi\Vert\\
&\leq 1+\delta|t|\\
\Psi(t=0)&=\psi
\endaligned\tag 3.41$$
where $\Vert\cdot\Vert$ can be either in $\Bbb L^2$ or analytic.  

Let $U=U(x, t)$ be the solution to the initial value problem:
$$\cases i U_t+\Delta U-[\delta (p+1)|u_\theta|^{2p}]U-[\delta p|u_\theta|^{2(p-1)}u_\theta^2]\bar U=0,\\
U(t=0)=\psi,\endcases\tag 3.42$$
where we have put back the suffix $\theta=\{\theta_j, j\in B_J\}$ to emphasize its dependence.
By definition
$$\align U(t)=&S(t)\psi\tag 3.43\\
=&S_\theta(t)\psi.\tag 3.44\endalign$$
Since $$\Psi(t)=U(t)+i\int_0^t S_\theta(t) S_\theta^{-1}(\tau)R(\tau)d\tau,\tag 3.45$$
we have 
$$S_\theta(t)\psi=U(t)=\Psi(t)-i\int_0^t S_\theta(t)S_\theta^{-1}(\tau)R(\tau)d\tau.\tag 3.46$$
Using $$S_{\theta+\omega't}(-t) S_\theta(t)=\Bbb I,\tag 3.47$$
it follows that 
$$\Vert S_\theta(t)\Vert\leq (1+\delta|t|)+|t|(\max_{\theta'\in\Bbb T^{|B_J|},|\tau|<|t|}\Vert S_{\theta'}(\tau)\Vert^2\delta^r
(1+|t|))\tag 3.48$$
implying that 
$$\Vert S_\theta(t)\Vert\leq 1+|t|\text{ for } |t|<\delta^{-r/3}.\tag 3.49$$
\hfill $\square$
\enddemo

\bigskip
\head{\bf 4. Almost global existence}\endhead
We now return to the Cauchy problem
$$
\cases iu_t =-\Delta u+\delta |u|^{2p}u,\qquad\, (p\geq 1, p\in\Bbb N \text{ \it arbitrary}),\,\,(4.1)\\
u(t=0)=u_0,\qquad\qquad\qquad\qquad\qquad\qquad\qquad\qquad\,\,(4.2)\endcases
$$
where $u_0\in\Bbb A$ satisfies $(\Bbb A1, 2)$ and prove the Theorem.
\demo{Proof of the Theorem}
The main idea is to find $v^{(0)}\in\Bbb A$ satisfying ($\Bbb A1, 2$) so that the corresponding quasi-periodic
solution $v$ constructed according to Proposition 2.1 satisfies:
$$
\cases iv_t +\Delta v-\delta |v|^{2p}v=\Cal O(\delta^r),\qquad\qquad\qquad\qquad\qquad\qquad(4.3)\\
v(t=0)-u_0=\Cal O(\delta^r),\qquad\qquad\qquad\qquad\qquad\qquad\qquad\,(4.4)\endcases
$$
the $\Cal O(\delta^r)$ is in the analytic norm in $(\Bbb A2)$ with $\beta''\in(0,\beta)$ (uniformly in $t$ for (4.3)). 
The quasi-periodic solution $v$ thus solves the Cauchy problem (4.1, 4.2) to $\Cal O(\delta^r)$. The control
of the flow linearized at $v$ will then prove the Theorem for $|t|<\delta^{-r}$.

(4.3, 4.4) follow from Proposition 2.1, Lemma 3.1 and the open mapping theorem as follows. Writing $u_0$
as $$u_0=u_1+u_2,\tag 4.5$$
where $u_1=\Cal O(1)$ is generic, $\{\hat u_1\}$ in $\Cal B''_\epsilon$ defined in Lemma 3.1, $u_2=\Cal O(\delta)$ and 
$$v^{(0)}=v_1+v_2\tag 4.6$$
as the initial approximation to the (approximate) quasi-periodic solution $v$, with $v_1$ generic, $\{\hat v_1\}$ in $\Cal B''_\epsilon$ defined in Lemma 3.1 and $v_2=\Cal O(\delta)$.

From (2.16) $$v_1=u_1,\tag 4.7$$
and hence we only need to determine $v_2$. Let
$$\beta=\Pi\hat u_2,\, \alpha=\Pi\hat v_2,\tag 4.8$$
where $$\cases (\Pi \hat w)(j)&= \hat w (j,-e_j), \quad |j|_\infty\leq\Cal O(|\log\delta|),\\
 &=0\qquad\qquad\qquad\qquad\text{otherwise}.\endcases$$
We take $e^{ij\cdot x}$ as the Fourier basis for $\Bbb T^d$.

Using Proposition 2.1, it suffces to look at the $P$-equations, since the $Q$-equations are satisfied and we are only interested
at $t=0$. (So $\omega t=0$.) We slightly modify the Newton scheme in (2.27) and let 
$$\aligned \Delta{\hat v}^{(1)}(\alpha)=&-\Pi {F_N'}^{-1}({\hat v}^{(0)}(\alpha))F({\hat v}^{(0)}(\alpha))\\
 \Delta{\hat v}^{(2)}(\alpha)=&-\Pi {F_N'}^{-1}({\hat v}^{(0)}+ \Delta{\hat v}^{(1)})(\alpha))F(({\hat v}^{(0)}+ \Delta{\hat v}^{(1)})(\alpha))\\
 &\vdots\endaligned\tag 4.9$$
and impose:
$$\sum_{k=1}^K\sum_n\Delta {\hat v}^{(k)}(\alpha)=\beta-\alpha,\tag 4.10$$
for some $K=K(r)$, which entails (4.3, 4.4), using also (4.7). 

Let $$\Cal F(\alpha)=\alpha+\sum_{k=1}^K\sum_n\Delta {\hat v}^{(k)}(\alpha).\tag 4.11$$
$\Cal F(\alpha)$ is smooth in $\alpha$ and moreover it follows from Lemma 3.1, (3.15, 3.25) that 
$$\Vert {\Cal F'}^{-1}(\alpha)\Vert<2\epsilon^{-1}\tag 4.12$$
for $\Vert\alpha\Vert=\Cal O(\delta)\ll\epsilon$. The open mapping theorem then gives that for any $\beta$, there is $\alpha$, so that 
(4.10) is satisfied, which in turn implies (4.3, 4.4). 

Decompose $u$ into 
$$u=v+w,\tag 4.13$$
with $v$ satisfying (4.3, 4.4). The remainder $w$ verifies 
$$\cases i w_t+\Delta w-[\delta (p+1)|v^{2p}]w-[\delta p|v|^{2(p-1)}v^2]\bar w+\delta f(w)=R=-\Cal O(\delta^r),\quad\quad(4.14)\\
w(t=0)=w_0=-\Cal O(\delta^r),\quad\quad\quad\quad\quad\quad\quad\quad\quad\quad\quad\qquad\qquad\quad\quad\quad\quad\qquad \,\,\,(4.15)\endcases$$
where $f(w)$ is at least quadratic, $f(w)=\Cal O(w^2)$ and the $\Cal O(\delta^r)$ are exactly as in (4.3, 4.4). 

The Duhamel formula then gives 
$$w(t)=S(t)w_0-i\int_0^t S(t)S^{-1}(\tau)[R(\tau)-\delta\Cal O(w^2(\tau))]d\tau,\tag 4.16$$
where $S(t)$ is the flow for the linear equation:
$$i g_t+\Delta g-[\delta (p+1)|v^{2p}]g-[\delta p|v|^{2(p-1)}v^2]\bar g=0.\tag 4.17$$

Let $w=\delta^{r'} w'$, $w_0=\delta^{r'}w_0'$ with some $0<r'<r$. Using (3.7), (4.16) gives 
$$\Vert w'(t)\Vert\leq (1+|t|)\Vert w_0'\Vert +(1+|t|^3)\delta^{r-r'}+(1+|t|^3)\delta^{r'} \Vert\sup w'\Vert^2.\tag 4.18$$
So $$\Vert w(t)\Vert <\Cal O(\delta^{r/2})\tag 4.19$$
for $|t|<\delta^{-r/10}$ from $\Vert w'(t)\Vert \leq 1$ by choosing $r'=r/2$. Combining with the estimates on $v$, we conclude the proof.
\hfill $\square$
\enddemo


\Refs\nofrills{References}
\widestnumber\key{CFKSA}

\ref
\key {\bf Ba}
\by D. Bambusi
\paper  Nekhoroshev theorem for small amplitude solutions in nonlinear Schr\"odinger equations
\jour Math. Z.
\vol 230
\pages 345-387
\yr 1999
\endref

\ref
\key {\bf BG}
\by D. Bambusi, B. Gr\'ebert
\paper  Birkhoff normal form for PDE's with tame modulus
\jour Duke Math. J. 
\vol 135
\pages 507-567
\yr 2006
\endref


\ref
\key {\bf Bo1}
\by J. Bourgain
\paper  Fourier transformation restriction phenomena for certain lattice subsets and applications to
nonlinear evolution equations, part I: Schr\"odinger equations
\jour Geom. and Func. Anal.
\vol 3
\pages 107-156
\yr 1993
\endref

\ref
\key {\bf Bo2}
\by J. Bourgain
\paper  Construction of approximative and almost periodic solutions of perturbed linear Schr\"odinger and wave equations
\jour Geom. and Func. Anal.
\vol 6
\pages 201-230
\yr 1996
\endref

\ref
\key {\bf Bo3}
\by J. Bourgain
\paper  On diffusion in high-dimensional Hamiltonian systems and PDE
\jour J. Anal. Math
\vol 80
\pages 1-35
\yr 2000
\endref

\ref
\key {\bf G}
\by M. Grassin
\paper  Global smooth solutions to Euler equations for a perfect gas
\jour Indiana Univ. Math. J. 
\vol 47
\pages 1397-1432
\yr 1998
\endref


\ref
\key {\bf Se}
\by D. Serre
\paper  Solution classique globales des \'equations d'Euler pour un fluide parfait compressible 
\jour Ann. Inst. Fourier
\vol 47
\pages 139-153
\yr 1997
\endref



\endRefs
\enddocument
\end